\documentclass{article}
\pdfoutput=1

\usepackage[narrow]{arxiv}

\usepackage[utf8]{inputenc} 
\usepackage[T1]{fontenc}    
\usepackage[hypertexnames=false]{hyperref}       
\usepackage{url}            
\usepackage{booktabs}       
\usepackage{amsfonts}       
\usepackage{nicefrac}       
\usepackage{microtype}      

\usepackage{lipsum}         
\usepackage{graphicx}
\usepackage{natbib}
\usepackage{doi}

\title{Optimal Projection-Free Adaptive SGD for Matrix Optimization}

\newif\ifuniqueAffiliation
\uniqueAffiliationtrue

\ifuniqueAffiliation 
\author{
  Dmitry Kovalev \\
  Yandex Research \\
  \texttt{dakovalev1@gmail.com}
}
\else
\usepackage{authblk}

\author[1]{%
  Dmitry Kovalev\thanks{\texttt{dakovalev1@gmail.com}}%
}
\affil[1]{Yandex Research}
\fi
\renewcommand{\shorttitle}{Optimal Projection-Free Adaptive SGD for Matrix Optimization}

\usepackage{algorithm,algpseudocodex}
\usepackage{nicematrix}

\usepackage[textsize=tiny,color=orange!20,bordercolor=orange,]{todonotes}
\presetkeys{todonotes}{inline}{}

\usepackage{xcolor}
\hypersetup{
  colorlinks=true,
  linkcolor=red!75!black,
  citecolor=blue!50!black
}

\usepackage{mathtools,amsthm,amssymb}
\allowdisplaybreaks[4]

\usepackage{enumitem}
\setlist{topsep=0pt,partopsep=0pt,itemsep=0pt}

\usepackage{cleveref}
\usepackage{crossreftools}
\crefname{problem}{problem}{problems}
\crefname{condition}{condition}{conditions}
\crefname{property}{property}{properties}
\creflabelformat{problem}{(#2#1#3)}

\usepackage{annotate}
\usepackage{mythm}
\usepackage{mymath}

\newcommand{\tmax}{{\textstyle\max}}
\newcommand{\tmin}{{\textstyle\min}}

\newcommand{\tg}{\tilde{g}}

\newcommand{\ox}{\overline{x}}
\newcommand{\hx}{\hat{x}}

\newcommand{\dist}[1]{\rho\brr{#1}}
\newcommand{\sqdist}[1]{\rho^2\brr{#1}}
\newcommand{\distd}[1]{\rho_*\brr{#1}}
\newcommand{\sqdistd}[1]{\rho_*^2\brr{#1}}

\newcommand{\sL}{\mathbb{L}}
\newcommand{\Spp}{\S_{++}}
\newcommand{\Sp}{\S_+}

\newcommand{\normo}{\gennorm{\mathrm{op}}{}}

\newcommand{\normt}{\gennorm{\mathrm{tr}}{}}

\newcommand{\N}{\mathbb{N}}
\newcommand{\bg}{\mathbb{B}}
\newcommand{\out}[1]{\mathop{\mathbf{out}}\brr{#1}}

\newcommand{\df}{\mathrm{d}}
\newcommand{\reg}{\mathrm{Reg}}

\begin{document}
\maketitle

\begin{abstract}
  Recently, \citet{jiang2026adaptive} developed Leon, a practical variant of One-sided Shampoo \citep{xie2025structured,an2025asgo} algorithm for online convex optimization, which does not require computing a costly quadratic projection at each iteration. Unfortunately, according to the existing analysis, Leon requires tuning an additional hyperparameter in its preconditioner and cannot achieve dimension-independent convergence guarantees for convex optimization problems beyond the bounded gradients assumption. In this paper, we resolve this issue by proving certain stability properties of Leon's preconditioner. Using our improved analysis, we show that tuning the extra hyperparameter can be avoided and, more importantly, develop {\em the first practical variant of One-sided Shampoo with Nesterov acceleration}, which does not require computing projections at each iteration. As a side contribution, we obtain improved dimension-independent rates in the non-smooth non-convex setting and develop a unified analysis of the proposed algorithm, which yields accelerated projection-free adaptive SGD with (block-)diagonal preconditioners.
\end{abstract}

\section{Introduction}

Stochastic gradient methods with diagonal (i.e., coordinate-wise) preconditioning, such as AdaGrad \citep{duchi2011adaptive}, Adam \citep{kingma2014adam}, and AdamW \citep{loshchilov2017decoupled}, have been considered the default choice for more than a decade in most deep learning tasks \citep{lecun2015deep}. However, recently, the focus of the research community has started to shift toward SGD algorithms with matrix preconditioning. The most notable examples of such algorithms are One-sided Shampoo \citep{xie2025structured} (also known as ASGO \citep{an2025asgo}), which incorporates gradient accumulation into the preconditioning matrix, and Muon \citep{jordan2024muon}, its accumulation-free alternative. The latter algorithm, Muon, has relatively cheap iterations and has already shown exceptional practical results in large-scale problems \citep{liu2025muon,team2025kimi}. On the other hand, One-sided Shampoo remains a highly promising algorithm. From a practical perspective, it shows good preliminary experimental results \citep{an2025asgo} while maintaining almost the same iteration cost as Muon. Moreover, from a theoretical perspective, One-sided Shampoo provably converges in the case of non-smooth objective functions \citep{an2025asgo} and is provably universal, i.e., it can adapt to different levels of H\"older smoothness with the same choice of hyperparameters \citep{kovalev2025sgd}, thanks to the gradient accumulation in its preconditioning matrix. However, there are still important open questions related to One-sided Shampoo that we aim to address in this paper.


\textbf{One-sided Shampoo vs AdaGrad.}
One-sided Shampoo was developed by \citet{xie2025structured} as a special instance of the adaptive meta-algorithm \citep{gupta2017unified}, which, in turn, is closely related to the original AdaGrad. The key feature of AdaGrad is that it considers the geometric properties of the optimization problems (such as distances between iterates or Lipschitz/smoothness constants of the objective function) with respect to the infinity vector norm, $\norm{}_\infty$, thanks to its diagonal preconditioning. This allows AdaGrad (and other diagonally preconditioned methods such as Adam/AdamW) to adapt to cases where the coordinates of the objective function have very different scales \citep{liu2024adagrad,jiang2024provable}. In contrast, One-sided Shampoo is used for solving optimization problems over a space of matrices and considers the geometric properties of the problems with respect to the spectral matrix norm, $\normo{}$, thanks to its more complex matrix preconditioner. This allows One-sided Shampoo to exploit more subtle geometric properties, such as low-rank gradients or approximately block-diagonal Hessians \citep{an2025asgo,xie2025structured}. Such properties are observed in deep neural networks \citep{zhao2021zero,yang2023spectral,cosson2023low,zhang2024transformers,zhang2024adam}, which may explain the recent success of the accumulation-free variant of One-sided Shampoo, Muon \citep{kovalev2025non,xie2025tale}.

\textbf{Projection-free Methods.} It is well known that both AdaGrad and One-sided Shampoo require the iterates to be bounded to provably converge \citep{gupta2017unified,xie2025structured,kovalev2025sgd}. To ensure this, a projection onto the appropriate (infinity or spectral) norm ball is typically performed at each iteration. However, while in AdaGrad this projection step is equivalent to a simple coordinate-wise clipping of the iterates, One-sided Shampoo requires solving a complex auxiliary problem of the form $\min_{\normo{X'} \leq 1} \<X'-X,Q(X'-X)>$, where $X,X' \in \R^{m\times n}$, and $Q \in \R^{m\times m}$ is a symmetric positive definite matrix. Recently, \citet{jiang2026adaptive} proposed Leon, a practical variant of One-sided Shampoo based on the Follow-the-Regularized-Leader (FTRL) algorithm \citep{abernethy2009competing,abernethy2016perturbation}, that completely eliminates the necessity for the costly projection. Unfortunately, this algorithm suffers from an important drawback: it requires tuning an additional scalar hyperparameter $\delta > 0$ (see \Cref{sec:recap}). According to the existing analysis, poor tuning of this parameter severely hurts the theoretical convergence guarantees and makes the resulting iteration complexity of Leon dimension-dependent. In addition, it requires the bounded gradients assumption, which significantly limits the range of applicability of Leon. Consequently, we arrive at the first main research question considered in this paper:
\begin{question}<q1>
  Is it possible to develop a practical projection-free variant of One-sided Shampoo that maintains dimension-independent convergence guarantees beyond the bounded gradients assumption and avoids tuning extra hyperparameters?
\end{question}

\textbf{Nesterov acceleration.}
It is well known that vanilla GD requires $\cO\brr{\epsilon^{-2/(1+\nu)}}$ iterations to find an $\epsilon$-accurate solution to a convex $\nu$-H\"older smooth minimization problem, where $\nu \in [0,1]$ \citep{nesterov2015universal}, and that this iteration complexity can be improved to $\cO\brr{\epsilon^{-2/(1+3\nu)}}$ with the help of Nesterov momentum \citep{nesterov1983method,nesterov2018lectures}. Similarly, the convergence of adaptive gradient methods can be accelerated with Nesterov momentum. In particular, an accelerated version of AdaGrad was developed by \citet{kovalev2025sgd}, and an accelerated variant of One-sided Shampoo was developed by \citet{xie2025tale}. Unfortunately, the accelerated One-sided Shampoo of \citet{xie2025tale} is not practical, as it still requires performing the complex projection step at each iteration. Consequently, we arrive at the second main research question considered in this paper:
\begin{question}<q2>
  Is it possible to develop a practical projection-free variant of One-sided Shampoo with Nesterov acceleration?
\end{question}

\textbf{Sumary of contributions.} In this paper, we provide positive answers to \Cref{q1,q2} and make the following contributions:
\begin{enumerate}[label*=\bf(\roman*)]
  \item We show that the preconditioning operator of Leon satisfies a certain {\em gradient stability property} in \Cref{lem:Psi3}. Using this property, we develop an improved analysis of the FTRL algorithm with Leon's preconditioning, which shows that the parameter $\delta$ can be set arbitrarily small. Therefore, we eliminate the necessity of tuning this parameter and obtain dimension-independent convergence guarantees for Leon in online and non-smooth non-convex optimization. Refer to \Cref{sec:FTRL} for details.

  \item We develop an accelerated version of the FTRL algorithm with Leon-like preconditioning, in which we incorporate the accumulation of squared gradient distances in a UniXGrad-like fashion \citep{kavis2019unixgrad,rodomanov2024universality}. In particular, we improve the analysis of UniXGrad to accommodate general preconditioners, whereas the previous analyses of UniXGrad were limited to scalar stepsizes. As a result, we develop the first accelerated and practical projection-free variant of One-sided Shampoo, which achieves the optimal $\cO\brr{\epsilon^{-2/(1+3\nu)}}$ iteration complexity for convex $\nu$-H\"older smooth optimization. Refer to \Cref{sec:convex} for details.

  \item As a side contribution, we provide a unified analysis of Leon and its accelerated version for optimization over general Euclidean spaces $\cX$. Therefore, our accelerated projection-free adaptive SGD algorithm is not limited to matrix preconditioners but can also incorporate a wide range of preconditioning operators, including diagonal preconditioners and scalar stepsizes.
\end{enumerate}

\textbf{Notation.}

In this paper, we use the following notation:
$\dim {\cX}$ is the dimension of the space $\cX$;
$\sL$ is the space of linear operators $\cX \to \cX$;
$\mI \in \sL$ and $\mO \in \sL$ denote the identity and the zero operators, respectively;
for arbitrary operator $\mA \in \sL$, $\mA^* \in \sL$ denotes its adjoint operator;
$\S \subset \sL$ is the space of self-adjoint linear operators, $\Spp, \Sp \subset \S$ are the spaces of positive definite and positive semi-definite self-adjoint operators, respectively; $\prec,\preceq,\succ,\succeq$ denote the standard L\"owner order on $\S$;
$\<,>$ and $\norm{}$ denote the standard inner product and the Euclidean norm on $\cX$ or $\sL$, depending on the context, in particular, $\<\mA,\mB> = \trace{\mA\mB^*}$ for $\mA,\mB \in \sL$, where $\trace{}$ is the standard trace functional on $\sL$ induced by the inner product $\<,>$ on $\cX$;
for arbitrary $\mH \in \Spp$, $\norm{}_\mH$ denotes the weighted Euclidean norm on $\cX$, i.e., $\sqn{x}_\mH = \<x,\mH x>$ for $x \in \cX$;
$\normo{}$ and $\normt{}$ denote the operator and trace norm on $\sL$, respectively, i.e., $\normo{\mA} = \max_{\norm{x}\leq 1}\norm{\mA x}$ and $\normt{\mA} = \trace{\sqrt{\mA\mA^*}}$ for all $\mA \in \sL$;
for arbitrary $y,z \in \cX$, by $z\<y,> \in \sL$ we denote the rank-1 linear operator $x\mapsto z\<y,x>$; $\out{z} \in \S$ denotes a symmetric rank-1 linear operator $z\<z,>$;
$\N$ and $\N_0$ denote the set of positive and non-negative integers.

\section{Preliminaries}\label{sec:pre}

Leon is a special instance of the FTRL algorithm for online convex optimization. Given a sequence of strongly convex closed regularization functions $\Psi_k(x)\colon \cX \to \R\cup\brf{+\infty}$ and a sequence of linear maps $x \mapsto \<g_k,x>$, where $g_k \in \cX$ are the past gradients, the update rule for FTRL is given as follows:
\begin{equation}\label{eq:FTRL}
  x_{k+1} = \argmin_{x \in \cX} \Psi_k(-x) + \<m_k,x>,\quad\text{where}\quad m_k = \msum_{i=0}^k g_i.
\end{equation}
From standard convex analysis, it follows that the Fenchel conjugate functions $\Psi_k^*(m)\colon \cX \to \cR$ are convex and smooth, and the iterations in \cref{eq:FTRL} can be written in a simplified form:
\begin{equation}\label{eq:FTRL2}
  x_{k+1} = -\nabla\Psi_k^*(m_k).
\end{equation}
Inspired by the choice of the regularization function proposed by \citet{jiang2026adaptive} for Leon, we will further use the following choice of $\Psi_k^*(m)$:
\begin{equation}\label{eq:Psi}
  \Psi_k^*(m) = \eta\normt{\sqrt{\proj[\cH]{\out{m}} + \mS_k}},
\end{equation}
where $\eta > 0$ is a positive parameter, $\mS_k \in \cH \cap \Spp$ is a positive definite operator, and $\cH \subset \S$ is a certain linear subspace of self-adjoint operators. The linear subspace $\cH$ in this definition is used for the purpose of unified analysis of the resulting algorithm. It is inspired by the unified analysis of the adaptive meta-algorithm by \citet{gupta2017unified}, where they restricted the choice of the preconditioners to the space $\cH$, e.g., diagonal and matrix preconditioners, as well as scalar stepsizes. Of course, to conduct the unified analysis, it is necessary to impose specific assumptions on the space $\cH$, e.g., to ensure that the operator under the square root in \cref{eq:Psi} is positive definite. There are two closely related sets of assumptions on $\cH$ proposed by \citet{xie2025structured,xie2025tale} and \citet{kovalev2025non,kovalev2025sgd}. We adhere to the latter approach with the following \Cref{ass:H}.
\begin{assumption}<ass:H>
  The linear subspace of operators $\cH \subset \S$ satisfies the following properties:
  \begin{enumerate}[label=\bf(\roman{*})]
    \item The space $\cH$ contains identity operator: $\mI \in \cH$.
    \item The space $\cH$ is closed under symmetric multiplications: $\mH\mF\mH \in \cH$ for all $\mH,\mF \in \cH$.
    \item The orthogonal projection onto $\cH$ is order preserving: $\proj[\cH]{\mH} \in \Spp$ for all $\mH \in \Spp$.
  \end{enumerate}
\end{assumption}
Following \citet[Lemma~1]{kovalev2025non}, we define the primal norm $\dist{} \colon \cX \to \R_+$ and the dual norm $\distd{}\colon \cX \to \R_+$, which are non-Euclidean in general, as follows:
\begin{equation}\label{eq:R}
  \dist{x} = \normo{\sqrt{\proj[\cH]{\out{x}}}}
  \quad\text{and}\quad
  \distd{x} = \normt{\sqrt{\proj[\cH]{\out{x}}}}.
\end{equation}
For instance, One-sided Shampoo and Leon assume $\cX = \R^{m\times n}$ and $\dist{}$ is the spectral matrix norm (up to constants), and diagonal AdaGrad assumes $\cX = \R^{d}$ and $\dist{}$ is the infinity vector norm. Refer to \citet{kovalev2025non,kovalev2025sgd,xie2025structured,xie2025tale} for more examples. These norms also play an important role in the description of the geometric properties of optimization problems that we will consider in this paper. In particular, we will consider the online convex optimization problems of the form
\begin{equation}\label{eq:regret}
  \text{minimize}\quad \reg_K = \brs*{\max_{x \in Q_\cR}\msum_{k=0}^K \<g_k,x_k - x>} \quad \text{s.t.}\quad x_0,\dots,x_k \in Q_\cR,
\end{equation}
where the decision set $Q_\cR = \brf{x \in \cX : \dist{x} \leq \cR}$ is the ball of radius $\cR > 0$ with respect to the norm $\dist{}$. We will also discuss the geometric properties of optimization problems related to \Cref{ass:H} in \Cref{sec:online,sec:problem}.

\section{Improved and Unified Analysis of Leon}\label{sec:FTRL}

\subsection{Unified Leon and Existing Convergence Results}\label{sec:recap}

We start with computing the gradient of the dual regularization function $\Psi_k^*(m)$ in the following \Cref{lem:Psi}, which generalizes the discussion of \citet[Section~4.2]{jiang2026adaptive} for the matrix optimization setting, $\cX = \R^{m\times n}$.

\begin{lemma}<lem:Psi>
  The function $\Psi_k^*(m)$ is convex and differentiable, and its gradient is given as follows:
  \begin{equation}
    \nabla\Psi_k^*(m) = \eta\brs*{\proj[\cH]{\out{m}} + \mS_k}^{-1/2}m.
  \end{equation}
\end{lemma}

Next, \citet[Theorem~4.2]{jiang2026adaptive} derived the feasibility, dominance, upper stability, and smoothness properties of the dual regularization function $\Psi_k^*(m)$ in the matrix optimization setting. We generalize these properties to the case of general $\cX$ under \Cref{ass:H} in the following \Cref{lem:Psi2}.

\begin{lemma}<lem:Psi2>
  Let $\mS_{k+1} \succeq \mS_k$. Then, for all $m,m' \in \cX$, the following inequalities hold:
  \begin{enumerate}[label=\bf(\roman{*})]
    \item Feasibility: $ \dist{\nabla \Psi_k^*(m)} \leq \eta$,
    \item Dominance: $\eta\distd{m} \leq\Psi_k^*(m)$,
    \item Upper stability: $0 \leq \Psi_{k+1}^*(m) - \Psi_k^*(m) \leq \eta\brs*{\normt{\sqrt{\mS_{k+1}}} - \normt{\sqrt{\mS_k}}}$,
    \item Smoothness: $\bg_{\Psi_k^*}(m;m') \leq \tfrac{1}{2}\eta\sqn{m-m'}_{\mS_k^{-1/2}}$.
  \end{enumerate}
\end{lemma}

Now, we are ready to present Leon for solving online problems of the form~\eqref{eq:regret}, which is formalized as \Cref{alg:1} and is unified according to the discussion in \Cref{sec:pre}. Note that the iterations of \Cref{alg:1} coincide with the FTRL iterations in \cref{eq:FTRL,eq:FTRL2}, where we use the squared gradient accumulation operator $\mS_k = \delta^2 \mI + \sum_{i=0}^{k}\proj[\cH]{\out{g_i}}$ in the definition of the dual regularization function $\Psi_k^*(m)$ in \cref{eq:Psi}. Also, note that the feasibility condition holds: $x_0,\dots,x_K \in Q_\cR$, due to \Cref{lem:Psi2}. We also restate the main result for this algorithm as \Cref{thm:baseline}. Note that \Cref{thm:baseline} is also presented in a generalized manner under \Cref{ass:H}, similar to the rest of this paper.

\begin{algorithm}[t]
  \caption{Generalized Leon for Online Optimization}
  \label{alg:1}
  \begin{algorithmic}[1]
    \State \bf{input:} $x_0 = 0 \in \cX$
    \State $\mS_{-1} = \delta^2 \mI, m_{-1} = 0$\label{line1:init}
    \For{$k=0,\dots,K$}
    \State observe gradient $g_k \in \cX$, $m_k = m_{k-1} + g_k$
    \State $\mS_k = \mS_{k-1} + \proj[\cH]{\out{g_k}}$\label{line1:S}
    \State $x_{k+1} = -\nabla \Psi_k^*(m_k)$\Comment{$\Psi_k^*(m)$ is defined in \cref{eq:Psi}}
    \EndFor
    \State {\bf output:} $(x_0,\dots,x_K) \in Q_\cR^{K+1}$
  \end{algorithmic}
\end{algorithm}

\begin{theorem}[Generalization of Theorem~3.2, \citet{jiang2026adaptive}]<thm:baseline>
  Let $\dist{g_k} \leq \delta$ for all $k \in \N_0$ and let $\eta = \cR$. Then, the output of \Cref{alg:1} satisfies the following inequality:
  \begin{equation}
    \reg_K \leq \cO\brs*{\delta\cR\dim{\cX} + \cR\distd{g_0} + \cR\normt{\sqrt{\mS_K}}}.
  \end{equation}
\end{theorem}

Now, we discuss the key problems with the result in \Cref{thm:baseline}, which are essential to solve in order to provide a positive answer to \Cref{q1}. First, \Cref{thm:baseline} requires the gradients to be uniformly bounded. This poses issues in the stochastic optimization case, where we typically assume the boundedness of the second moments of the gradients, which does not imply that the gradients are almost surely bounded. Second, the presence of the term $\delta\cR \dim{\cX}$ in the regret upper bound prevents us from obtaining dimension-independent complexities for solving smooth convex problems and from accelerating the convergence of the algorithm using Nesterov momentum.

\subsection{Improved Analysis of Leon and FTRL}

First, for the FTRL algorithm in \cref{eq:FTRL2}, instead of bounding the standard regret defined in \cref{eq:regret}, we will bound the regret-like quantity $\reg_K^+$ defined as follows:
\begin{equation}\label{eq:regret2}
  \reg_K^+ = \max_{x \in Q_\cR}\brs*{\msum_{k=0}^K \<g_k,x_{k+1} - x>}.
\end{equation}

We obtain the key identity for the regret-like quantity $\reg_K^+$ in the following \Cref{lem:regret+}. In contrast to the analogous identity in Lemma~3.1 of \citet{jiang2026adaptive}, the inequality in \Cref{lem:regret+} contains negative Bregman divergence terms $-\bg_{\Psi_k^*}(m_{k-1};m_k)$, which we will be able to utilize by applying the smoothness property in \Cref{lem:Psi2}.

\begin{lemma}<lem:regret+>
  Let $\eta = \cR$ and $\delta > 0$. Then, the iterations in \cref{eq:FTRL2} satisfy the following relation:
  \begin{equation}
    \begin{aligned}
      \reg_K^+ &= \Psi_0^*(0) + \eta \distd{m_K} - \Psi_K^*(m_K)
      \\&
      + \msum_{k=1}^K \brs*{\Psi_k^*(m_{k-1}) - \Psi_{k-1}^*(m_{k-1})}
      - \msum_{k=0}^K \bg_{\Psi_k^*}(m_{k-1};m_k),
    \end{aligned}
  \end{equation}
\end{lemma}

Next, we establish the key property of the dual regularization function $\Psi_k^*(m)$, which we call gradient stability, in the following \Cref{lem:Psi3}.

\begin{lemma}<lem:Psi3>
  Let $\mS_{k+1} \succeq \mS_k$. Then, for all $m \in \cX$, the following inequality holds:
  \begin{equation}
    \sqn{\nabla \Psi_{k+1}^*(m) - \nabla \Psi_k^*(m)}_{\mS_{k+1}^{1/2}}
    \leq
    \eta^2\brs*{\normt{\sqrt{\mS_{k+1}}} - \normt{\sqrt{\mS_k}}}.
  \end{equation}
\end{lemma}

Now, using the previously established properties in \Cref{lem:Psi2}, the new property in \Cref{lem:Psi3}, and the negative Bregman divergence terms $-\bg_{\Psi_k^*}(m_{k-1};m_k)$, we are ready to establish the main upper bound on the regret-like quantity $\reg_K^+$ in the following \Cref{thm:FTRL}.

\begin{theorem}<thm:FTRL>
  Let $\eta = \cR$ and $\delta > 0$. Let $\mS_{k+1} \succeq \mS_k$ for all $k \in \N_0$. Let $g_0,\dots,g_K \in \cX$ and $x_0 = 0$. Let $x_1,\dots,x_{K+1}$ be generated according to \cref{eq:FTRL2,eq:Psi}. Then, it holds that
  \begin{equation}
    \reg_K^+ \leq \Psi_0^*(0) + \tfrac{3}{2}\eta\normt{\sqrt{\mS_K}}
    - \msum_{k=0}^K \tmfrac{1}{4\eta}\sqn{x_{k+1} - x_k}_{\mS_k^{1/2}}.
  \end{equation}
\end{theorem}

In the following section, we show how to use the upper-bound in \Cref{thm:FTRL} for obtaining improved guarantees in online convex and non-smooth non-convex optimization.

\subsection{Application to Online and Non-smooth Non-convex Optimization}\label{sec:online}

\textbf{Online optimization.}
We start with the general improved result for \Cref{alg:1} for solving online convex optimization, as stated in the following \Cref{thm:online}.

\begin{theorem}<thm:online>
  Let $\eta = \cR$ and $\delta > 0$. Then, the output of \Cref{alg:1} satisfies the inequality
  \begin{equation}
    \reg_K
    \leq
    \cO\brs*{\delta\cR\dim{\cX} + \cR\distd{g_0} + \cR\normt{\sqrt{\mS_K}}}.
  \end{equation}
\end{theorem}

The upper bound on the regret in the improved \Cref{thm:online} coincides with the inequality in \Cref{thm:baseline} by \citet{jiang2026adaptive}. However, we significantly simplify the assumptions required to achieve this inequality. First, \Cref{thm:online} does not require the gradients $g_k$ to be uniformly bounded. Second, we do not have to tune the parameter $\delta > 0$, as we can simply choose an arbitrarily small value for $\delta$, thus obtaining dimension-independent guarantees in \Cref{thm:online}.

\textbf{Non-smooth non-convex optimization.}
Now, we show, how to apply the improved result for solving optimization problems of the form
\begin{equation}
  \min_{x \in \cX} f(x),
\end{equation}
where $f(x)\colon \cX \to \cR$ is a differentiable objective function. Our goal is to find an $(\gamma,\epsilon)$-stationary point $\hx \in \cX$, i.e., satisfying the following \Cref{def:stationary}.

\begin{definition}<def:stationary>
  A point $\hx \in \cX$ is called $(\gamma,\epsilon)$-stationary if there exists a distribution $\cP$ over $\cX$ such that \textbf{(i)} $\E[z \sim \cP]{z} = \hx$, \textbf{(ii)} $\distd{\E[z \sim \cP]{\nabla f(z)}} \leq \epsilon$, and \textbf{(iii)} $\E[z \sim \cP]{\dist{z - \hx}} \leq \gamma$.
\end{definition}

We are naturally interested in solving the problem using a stochastic algorithm. Hence, we assume the existence of an appropriate stochastic gradient estimator in the following \Cref{ass:grad1}.

\begin{assumption}<ass:grad1>
  There exists a stochastic estimator $\nabla_\xi f(x)$ of the gradient $\nabla f(x)$, where $\xi \sim \cD$ is a random variable sampled from the distribution $\cD$, which satisfies the following properties:
  \begin{enumerate}[label=\bf(\roman*)]
    \item Unbiased estimator: $\E[\xi \sim \cD]{\nabla_\xi f(x)} = \nabla f(x)$.
    \item 2nd moment bound: $\E[\xi \sim \cD]{\sqn{\nabla_\xi f(x)}_{\mB^{-1}}} \leq \cG^2$, where $\mB \in \Spp\cap \cH$, $\normt{\mB} = 1$, $\cG > 0$.
  \end{enumerate}
\end{assumption}

Now, we can state the complexity result for \Cref{alg:1} for solving non-smooth non-convex minimization problems. It is obtained with the help of the O2NC framework \citep{zhang2024random}. The proof and the resulting algorithm are almost identical to those described by \citet[Section~6]{jiang2026adaptive}; hence, we will not describe them in this paper. Note that, similar to the online optimization case, we obtain improved dimension-independent guarantees in \Cref{thm:nonconvex} by choosing a very small value of $\delta > 0$.

\begin{theorem}<thm:nonconvex>
  \Cref{alg:1} with the O2NC framework requires the following number of iterations to find an $(\gamma,\epsilon)$-stationary point, where $\Delta = f(x_0) - \inf_x f(x)$ is the initial function gap:
  \begin{equation}
    \cO\brs*{\Delta\cG^2 \cdot  \gamma^{-1}\epsilon^{-3}+ \cG^2\cdot\epsilon^{-2} + \delta \dim{\cX}\cdot\epsilon^{-1}}.
  \end{equation}
\end{theorem}

\textbf{On the geometry of \Cref{ass:grad1}.}
Note that \Cref{ass:grad1} is given in terms of the weighted Euclidean norm $\norm{}_{\mB^{-1}}$. In the case of matrix optimization, $\cX = \R^{m\times n}$, it is convenient to use it for modeling low-rank gradients, which are often observed in deep neural networks. Indeed, \Cref{ass:grad1} is implied by the inequality $\E[\xi \sim \cD]{G_\xi G_\xi^\top}\preceq C^2$, where $G_\xi \in \R^{m\times n}$ is a stochastic gradient, and $C \in \R^{m\times m}$ is a symmetric positive definite matrix with low stable rank.\footnote{By the stable rank of a symmetric positive definite matrix $C$, we mean $\trace{C} / \lambda_{\max}(C)$.} More details are available in the work of \citet[Section~4]{an2025asgo}. In addition, in the general case, \Cref{ass:grad1} is also closely related to its non-Euclidean analogue: it implies the inequality $\E[\xi \sim \cD]{\sqdistd{\nabla_\xi f(x)}} \leq \cG^2$. Note that optimization under non-Euclidean geometry assumptions has attracted a lot of interest in the research community \citep{bernstein2024old,kovalev2025understanding,pethick2025training,kovalev2025non}.

\section{Generalized Leon with Nesterov Acceleration for Convex Problems}\label{sec:convex}

\subsection{Problem and Assumptions}\label{sec:problem}

In this section, we focus on solving convex optimization problems of the form
\begin{equation}\label[problem]{eq:main}
  f^* = \min_{x \in Q_\cR} f(x).
\end{equation}
We start by imposing the formal \Cref{ass:f,ass:grad2} on the objective function $f(x)$. These assumptions can be seen as a generalization of the previously discussed \Cref{ass:grad1} to the case of stochastic H\"older smooth convex optimization. Consequently, these assumptions can also be used to model low-rank gradients in the matrix optimization case, as discussed in \Cref{sec:online}. It is also important to highlight that, in the matrix case, the constraint $x \in Q_\cR$ implies that the solution is likely to be high-rank, as discussed by \citet{chen2025muon}, which can be beneficial, as discussed by \citet[Section~5]{an2025asgo}. Additionally, in the general case, \Cref{ass:f} implies H\"older smoothness of the objective function $f(x)$ with respect to the non-Euclidean norm $\dist{}$, which makes it relevant to non-Euclidean optimization \citep{bernstein2024old,kovalev2025understanding,pethick2025training,kovalev2025non}.

\begin{assumption}<ass:f>
  Function $f(x)$ is convex and $(\cL,\nu)$-H\"older smooth with respect to the norm $\norm{}_\mB$, where $\cL > 0$, $\nu \in [0,1]$, $\mB \in \Spp\cap \cH$ and $\normt{\mB} = 1$.
\end{assumption}

\begin{assumption}<ass:grad2>
  There exists a stochastic estimator $\nabla_\xi f(x) = n_\xi(x) + \nabla f(x)$ of the (sub)gradient $\nabla f(x) \in \partial f(x)$, where $n_\xi(x)$ is the noise and $\xi \sim \cD$ is a random variable sampled from the distribution $\cD$. The noise $n_\xi(x)$ satisfies the following properties:
  \begin{enumerate}[label=\bf(\roman*)]
    \item Zero mean: $\E[\xi \sim \cD]{n_\xi(x)} = 0$.
    \item Variance bound: $\E[\xi \sim \cD]{\sqn{n_\xi(x)}_{\mB^{-1}}} \leq \sigma^2$, where $\mB \in \Spp\cap \cH$, $\normt{\mB} = 1$, and $\sigma > 0$.
  \end{enumerate}
\end{assumption}

\begin{algorithm}[t]
  \caption{Generalized Leon with Nesterov Acceleration}
  \label{alg:2}
  \begin{algorithmic}[1]
    \State {\bf input:} $x_0 = \ox_0 = 0 \in Q_\cR$
    \State $\mS_0 = \delta^2 \mI$, $m_{-1} = 0$
    \label{line2:init}
    \For{$k=0,\dots,K$}
    \State sample $\xi_k, \xi_{k+1}' \sim \cD$
    \State $g_k = \nabla_{\xi_k} f_k(x_k)$, $m_k = m_{k-1} + g_k$ \label{line2:g1}
    \Comment{$f_k(x)$ is defined in \cref{eq:f_k}}
    \State $x_{k+1} = -\nabla \Psi_k^*(m_k)$\Comment{$\Psi_k^*(m)$ is defined in \cref{eq:Psi}}
    \State $\ox_{k+1} = x_{k+1}/\alpha_k + (1-1/\alpha_k)\ox_k$\label{line2:ox}
    \State $\tg_{k+1} = \nabla_{\xi'_{k+1}} f_k(x_{k+1})$ \label{line2:g2}
    \Comment{$f_k(x)$ is defined in \cref{eq:f_k}}
    \State $\mS_{k+1} = \mS_{k} + \proj[\cH]{\out{\tg_{k+1} - g_k}}$ \label{line2:S}
    \EndFor
    \State {\bf output:} $\ox_{K+1} \in Q_\cR$
  \end{algorithmic}
\end{algorithm}

\subsection{Optimal Projection-Free Adaptive SGD with Preconditioning}

Now, we are ready to present the accelerated version of Leon with Nesterov acceleration, formalized as \Cref{alg:2}. In order to incorporate Nesterov acceleration into the algorithm, we follow the approach of \citet{kovalev2024linear}. That is, at each iteration $k$, we apply the FTRL step in \cref{eq:FTRL2} to the modified objective function $f_k(x)\colon \cX \to \R$, which is defined as follows:
\begin{equation}\label{eq:f_k}
  f_k(x) = \alpha_k^2 f(x/\alpha_k + (1-1/\alpha_k)\ox_k),
  \quad\text{where}\;\; \alpha_k \geq 1,\;\;\ox_k \in \cX,
\end{equation}
Consequently, using our previously obtained \Cref{thm:FTRL}, we obtain the following \Cref{lem:function_gap}.

\begin{lemma}<lem:function_gap>
  Let $\eta = \cR$ and $\delta > 0$. Then, for all $x \in Q_\cR$, the following inequality holds:
  \begin{equation}
    \E*{\msum_{k=0}^K\bg_{f_k}(x_k;x_{k+1}) + f_k(x_{k+1}) - f_k(x)} \leq \E*{\delta\cR \dim{\cX} + 4\cR\normt{\sqrt{\mS_{K+1}}}}.
  \end{equation}
\end{lemma}

The reference point $\ox_k$ is updated using \cref{line2:ox}, which allows us to obtain the following \Cref{lem:acceleration}.

\begin{lemma}<lem:acceleration>
  Let $\alpha_0 = 1$ and $\alpha_k^2 \leq \alpha_k + \alpha_{k-1}^2$ for all $k \in \N$. Let $x^* \in Q_\cR$ be a solution to \cref{eq:main}. Then, the following inequality holds:
  \begin{equation}
    \alpha_K^2 \brs*{f(\ox_{K+1}) - f^*} \leq \msum_{k=0}^K\brs*{f_k(x_{k+1}) - f_k(x^*)}.
  \end{equation}
\end{lemma}

Furthermore, we make another important modification to the algorithm. It is related to the fact that the convergence of standard AdaGrad-like algorithms relies on bounding the squared gradients using the objective function optimality gap \citep[Corollary~1]{orabona2023normalized}:
\begin{equation}
  \norm{\nabla f(x_k)}_{\mB^{-1}}^{1+1/\nu} \leq \tmfrac{1+\nu}{\nu} \cdot \cL^{1/\nu} \cdot \brs{f(x_k) - \tmin_{x \in \cX} f(x)}.
\end{equation}
However, the minimizer of the objective function $f(x)$ does not coincide with the minimizer of the modified function $f_k(x)$, even in the unconstrained optimization case, which prevents us from using this inequality in an accelerated algorithm. To tackle this issue, we accumulate squared gradient differences in the operator $\mS_k$ in \cref{line2:S}, instead of merely accumulating squared gradients, similar to the UniXGrad algorithm. In addition, we can bound the squared gradient differences with the Bregman divergences as follows:
\begin{equation}\label{eq:grad_diff}
  \norm{\nabla f_k(x_{k+1}) - \nabla f_k(x_k)}_{\mB^{-1}}^{1+1/\nu} \leq \tmfrac{1+\nu}{\nu} \cdot \brs*{\cL\alpha_k^{1-\nu}}^{1/\nu} \cdot \bg_{f_k}(x_k;x_{k+1}).
\end{equation}
Consequently, we can bound the squared gradient differences by utilizing the negative Bregman divergence term in \Cref{lem:function_gap}. This improves upon the previous analyses by \citet{kavis2019unixgrad,rodomanov2024universality}, who did not properly utilize this term and achieved convergence only for scalar stepsizes. In accordance with this discussion, we obtain the following \Cref{lem:S}.

\begin{lemma}<lem:S>
  Let $\alpha_k \geq \alpha_{k-1}$ for all $k \in \N$. Then, for all $c \geq 1$, the following inequality holds:
  \begin{equation}
    \begin{aligned}
      \E*{\normt{\sqrt{\mS_{K+1}}}}
      \leq
      \delta\dim{\cX}
      &+
      \E*{\tmfrac{1}{c\eta}\msum_{k=0}^K\bg_{f_k}(x_k;x_{k+1})}
      \\&+
      2\sigma\brs*{[K+1]\alpha_K^2}^{\frac{1}{2}}
      +
      c\cL\eta^{\nu} \brs*{[K+1]\alpha_K^2}^{\frac{1-\nu}{2}}.
    \end{aligned}
  \end{equation}
\end{lemma}

Now, it remains to combine \Cref{lem:function_gap,,lem:acceleration,,lem:S} in order to derive the final convergence result for \Cref{alg:2}, as stated in the following \Cref{thm:convex}.

\begin{theorem}<thm:convex>
  Let $\eta = \cR$ and $\delta > 0$, and let $\alpha_k = 1 + k/2$. Then, the output of \Cref{alg:2} satisfies the following inequality:
  \begin{equation}
    \E{f(\ox_{K+1}) - f(x^*)} \leq
    \cO\brs*{
      \tmfrac{\cL\cR^{1+\nu}}{[K+1]^{(3\nu+1)/2}}
      +
      \tmfrac{\sigma\cR}{[K+1]^{1/2}}
      +
      \tmfrac{\delta\cR\dim{\cX}}{[K+1]^2}
    }.
  \end{equation}
  Hence, to reach the precision $\E{f(\ox_{K+1}) - f^*} \leq \epsilon$, it is sufficient to perform the following number of iterations of \Cref{alg:2}:
  \begin{equation}
    K + 1 = \cO\brs*{\brs*{\tmfrac{\cL\cR^{1+\nu}}{\epsilon}}^{\frac{2}{1+3\nu}} + \brs*{\tmfrac{\sigma\cR}{\epsilon}}^2 + \brs*{\tmfrac{\delta\cR \dim{\cX}}{\epsilon}}^{\frac{1}{2}}}.
  \end{equation}
\end{theorem}

\Cref{thm:convex} implies the optimal complexity $\cO\brr{\epsilon^{-2/(3\nu+1)}}$ for non-stochastic $\nu$-H\"older smooth convex problems and the optimal complexity $\cO\brr{\epsilon^{-2}}$ for stochastic convex optimization problems \citep{nemirovskij1983problem}. In addition, the complexity result in \Cref{thm:convex} is dimension-independent by choosing a very small value of $\delta > 0$. Moreover, it involves the constants $\sigma, \cL > 0$ from \Cref{ass:f,ass:grad2}, and the radius $\cR > 0$ of the constraint set $Q_\cR$ with respect to the non-Euclidean norm $\dist{}$. Therefore, in the matrix optimization case, our \Cref{alg:2} can exploit low-rank gradients and a high-rank solution to outperform its Euclidean counterparts, just like the original non-accelerated One-sided Shampoo \citep[Section~5]{an2025asgo}.

\newpage

\cite*{}

\bibliographystyle{unsrtnat}
\bibliography{references}

\newpage
\appendix

\section[Proofs for \crtCref{sec:FTRL}]{Proofs for \Cref{sec:FTRL}}

\proofsubsection{lem:Psi}

We compute the differential $\df \Psi_k^*(m)$ as follows:
\begin{align*}
  \df \Psi_k^*(m)
  &\at{uses the differentiation rule for trace functions}=
  \tfrac{1}{2}\eta \<\brs{\proj[\cH]{\out{m}} + \mS_k}^{-1/2}, \df \brs{\proj[\cH]{\out{m}} + \mS_k}>
  \\&\at{uses the linearity of $\proj[\cH]{}$ and the bilinearity of $\out{}$}=
  \tfrac{1}{2}\eta \<\brs{\proj[\cH]{\out{m}} + \mS_k}^{-1/2}, \proj[\cH]{\df m\<m,> + m\<\df m,>}>
  \\&\at{uses the fact that $\proj[\cH]{\out{m}} + \mS_k \in \cH$, \Cref{ass:H}, and the properties of the projection}=
  \tfrac{1}{2}\eta \<\brs{\proj[\cH]{\out{m}} + \mS_k}^{-1/2}, \brs{\df m\<m,> + m\<\df m,>}>
  \\&= \<\eta\brs{\proj[\cH]{\out{m}} + \mS_k}^{-1/2}m, \df m>,
\end{align*}
where \annotate.\qed

\proofsubsection{lem:Psi2}

\textbf{Feasibility.}
\Cref{lem:Psi} implies $\nabla\Psi_k^*(m) = \eta \mA_k^{-1/2}m$, where $\mA_k = \proj[\cH]{\out{m}} + \mS_k$. Hence, we can upper-bound $\sqdist{\nabla\Psi_k^*(m)}$ as follows:
\begin{align*}
  \sqdist{\nabla\Psi_k^*(m)}
  &\at{uses the definition in \Cref{eq:R}}=
  \normo{\proj[\cH]{\out{\nabla\Psi_k^*(m)}}}
  =
  \tmax_{\norm{u} \leq 1} \<\proj[\cH]{\out{\nabla\Psi_k^*(m)}}, \out{u}>
  \\&\at{uses the properties of the projection}=
  \tmax_{\norm{u} \leq 1} \<\out{\nabla\Psi_k^*(m)}, \proj[\cH]{\out{u}}>
  \\&\at{uses the expression for $\nabla\Psi_k^*(m)$ above}=
  \eta^2\tmax_{\norm{u} \leq 1} \<\mA_k^{-1/2}\out{m}\mA_k^{-1/2}, \proj[\cH]{\out{u}}>
  \\&=
  \eta^2\tmax_{\norm{u} \leq 1} \<\out{m}, \mA_k^{-1/2}\proj[\cH]{\out{u}}\mA_k^{-1/2}>
  \\&\at{uses the fact that $\mA_k^{-1/2} \in \cH$, \Cref{ass:H} and the properties of the projection}=
  \eta^2\tmax_{\norm{u} \leq 1} \<\proj[\cH]{\out{m}}, \mA_k^{-1/2}\proj[\cH]{\out{u}}\mA_k^{-1/2}>
  \\&\at{uses the fact that $\proj[\cH]{\out{m}} \preceq \mA_k$ and the fact that $\mA_k^{-1/2}\proj[\cH]{\out{u}}\mA_k^{-1/2} \succeq \mO$ due to \Cref{ass:H}}\leq
  \eta^2\tmax_{\norm{u} \leq 1} \<\mA_k, \mA_k^{-1/2}\proj[\cH]{\out{u}}\mA_k^{-1/2}>
  \\&=
  \eta^2\tmax_{\norm{u} \leq 1} \<\mI, \proj[\cH]{\out{u}}>
  \\&\at{uses the fact that $\mI \in \cH$ due to \Cref{ass:H} and the properties of the projection}=
  \eta^2\tmax_{\norm{u} \leq 1} \<\mI, \out{u}>
  =
  \eta^2\tmax_{\norm{u} \leq 1} \sqn{u}
  =
  \eta^2,
\end{align*}
where \annotate.

\textbf{Dominance.}
We can lower-bound $\Psi_k^*(m)$ as follows:
\begin{align*}
  \Psi_k^*(m)
  &=
  \eta\<\mA_k^{1/2}, \mI>
  \at{uses the monotonicity of the trace square root}\geq
  \eta\<\sqrt{\proj[\cH]{\out{m}}}, \mI>
  \at{uses the definition in \cref{eq:R}}=
  \eta\distd{m},
\end{align*}
where \annotate.

\textbf{Upper stability.} The inequality $\Psi_{k+1}^*(m) \geq \Psi_k^*(m)$ is implied by the monotonicty of the trace square root. Furthermore, $\Psi_{k+1}^*(m)$ can be upper-bounded as follows:
\begin{align*}
  \Psi_{k+1}^*(m)
  &\at{uses the definition in \cref{eq:Psi}}=
  \eta\<\mI, \mA_{k+1}^{1/2} - \mA_k^{1/2}>
  \\&\at{use the differentiability properties of the trace square root}\leq
  \int_0^1 \tfrac{1}{2}\eta\<\brs{\mA_k + t(\mA_{k+1} - \mA_k)}^{-1/2},\mA_{k+1} - \mA_k>\df t
  \\&=
  \int_0^1 \tfrac{1}{2}\eta\<\brs{\mA_k + t(\mS_{k+1} - \mS_k)}^{-1/2},\mS_{k+1} - \mS_k>\df t
  \\&\at{uses the assumption $\mS_{k+1} \succeq \mS_k$, the fact that $\mA_k \succeq \mS_k$, and the L\"owner Heinz theorem \citep[Theorem~2.6]{carlen2010trace}}=
  \int_0^1 \tfrac{1}{2}\eta\<\brs{\mS_k + t(\mS_{k+1} - \mS_k)}^{-1/2},\mS_{k+1} - \mS_k>\df t
  \\&\at{use the differentiability properties of the trace square root}\leq
  \eta\<\mI, \mS_{k+1}^{1/2} - \mS_k^{1/2}>
  =
  \eta\brs*{\normt{\sqrt{\mS_{k+1}}} - \normt{\sqrt{\mS_k}}},
\end{align*}
where \annotate.

\textbf{Smoothness.}
It is sufficient to upper-bound the second-order differential $\df^2 \Psi_k^*(m)$:
\begin{align*}
  \df^2 \Psi_k^*(m)
  &\at{uses the definition in \cref{eq:Psi}}=
  \df^2\brs{\eta\<\mI, \mA_k^{1/2}>}
  \\&\at{uses the differentiability properties of the trace square root}=
  \tfrac{1}{2}\eta\df\brs{\<\mA_k^{-1/2},\df \mA_k>}
  =
  \tfrac{1}{2}\eta\<\mA_k^{-1/2},\df^2 \mA_k>
  +\tfrac{1}{2}\eta\<\df\mA_k^{-1/2},\df \mA_k>
  \\&\at{uses the definition of $\mA_k$ and the linearity of the projection}=
  \tfrac{1}{2}\eta\<\mA_k^{-1/2},\proj[\cH]{2\out{\df m}}>
  +\tfrac{1}{2}\eta\<\df\mA_k^{-1/2},\df \mA_k>
  \\&\at{uses the fact that $\mA_k^{-1/2} \in \cH$ due to \Cref{ass:H}}=
  \eta\<\mA_k^{-1/2},\out{\df m}>
  +\tfrac{1}{2}\eta\<\df\mA_k^{-1/2},\df \mA_k>
  \\&=
  \eta\sqn{\df m}_{\mA_k^{-1/2}}
  +\tfrac{1}{2}\eta\<\df\mA_k^{-1/2},\df \mA_k>
  \\&\at{uses the fact that $\tfrac{1}{2}\<\df\mA_k^{-1/2},\df \mA_k>$ is the second-order differential of the concave function $\mA \mapsto \<\mA^{1/2},\mI>$}\leq
  \eta\sqn{\df m}_{\mA_k^{-1/2}}
  \at{uses the fact that $\mS_k \preceq \mA_k$ and the L\"owner-Heinz theorem \citep[Theorem~2.6]{carlen2010trace}}\leq
  \eta\sqn{\df m}_{\mS_k^{-1/2}},
\end{align*}
where \annotate.\qed

\proofsubsection{lem:regret+}

We can rewrite $\reg_K^+$ as follows:
\begin{align*}
  \reg_K^+
  &\at{uses the definition in \cref{eq:regret2} duality between the norms defined in \cref{eq:R}}=
  \msum_{k=0}^K \<g_k,x_{k+1}> + \eta \distd{m_K}
  \\&\at{uses \cref{eq:FTRL2}}=
  \msum_{k=0}^K \<g_k,-\nabla \Psi_k^*(m_k)> + \eta \distd{m_K}
  \\&\at{uses \cref{eq:FTRL}}=
  \msum_{k=0}^K \<m_{k-1} - m_k,\nabla \Psi_k^*(m_k)> + \eta \distd{m_K}
  \\&=
  \msum_{k=0}^K \brs*{\Psi_k^*(m_{k-1}) - \Psi_k^*(m_k) -\bg_{\Psi_k^*}(m_{k-1};m_k)} + \eta \distd{m_K}
  \\&=
  \Psi_0^*(m_{-1}) - \Psi_0^*(m_0)
  +  \msum_{k=1}^K \brs*{\Psi_k^*(m_{k-1}) - \Psi_{k-1}^*(m_{k-1}) + \Psi_{k-1}^*(m_{k-1}) - \Psi_k^*(m_k)}
  \\&
  +\eta \distd{m_K} - \msum_{k=0}^K \bg_{\Psi_k^*}(m_{k-1};m_k)
  \\&=
  \Psi_0^*(0) + \eta \distd{m_K} - \Psi_K^*(m_K)
  + \msum_{k=1}^K \brs*{\Psi_k^*(m_{k-1}) - \Psi_{k-1}^*(m_{k-1})}
  \\&
  - \msum_{k=0}^K \bg_{\Psi_k^*}(m_{k-1};m_k),
\end{align*}
where \annotate.\qed

\proofsubsection{lem:Psi3}
Let $\mA_k = \proj[\cH]{\out{m}}+\mS_k$. Then, we can obtain the following:
\begin{align*}
  \mi{3}\sqn{\nabla \Psi_{k+1}^*(m) - \nabla \Psi_k^*(m)}_{\mS_{k+1}^{1/2}}
  \\&\at{uses \Cref{lem:Psi}}=
  \eta^2\sqn{(\mA_{k+1}^{-1/2} - \mA_k^{-1/2})m}_{\mS_{k+1}^{1/2}}
  \\&\at{uses the fact that $\mA_{k+1} \succeq \mS_{k+1}$ and the L\"owner-Heinz theorem \citep[Theorem~2.6]{carlen2010trace}}\leq
  \eta^2\sqn{(\mA_{k+1}^{-1/2} - \mA_k^{-1/2})m}_{\mA_{k+1}^{1/2}}
  \\&=
  \eta^2\<(\mA_{k+1}^{-1/2} - \mA_k^{-1/2})\mA_{k+1}^{1/2}(\mA_{k+1}^{-1/2} - \mA_k^{-1/2}),\out{m}>
  \\&\at{uses the properties of the projection and the inclusion $(\mA_{k+1}^{-1/2} - \mA_k^{-1/2})\mA_{k+1}^{1/2}(\mA_{k+1}^{-1/2} - \mA_k^{-1/2}) \in \cH$ due to the definition of $\mA_k$, \Cref{ass:H}, and the inclusion $\mS_k \in \cH$}=
  \eta^2\<(\mA_{k+1}^{-1/2} - \mA_k^{-1/2})\mA_{k+1}^{1/2}(\mA_{k+1}^{-1/2} - \mA_k^{-1/2}),\proj[\cH]{\out{m}}>
  \\&\at{uses the fact that $(\mA_{k+1}^{-1/2} - \mA_k^{-1/2})\mA_{k+1}^{1/2}(\mA_{k+1}^{-1/2} - \mA_k^{-1/2}) \succeq \mO$  and  $\proj[\cH]{\out{m}} \preceq\mA_k$}\leq
  \eta^2\<(\mA_{k+1}^{-1/2} - \mA_k^{-1/2})\mA_{k+1}^{1/2}(\mA_{k+1}^{-1/2} - \mA_k^{-1/2}),\mA_k>
  \\&=
  \eta^2\<\mA_{k+1}^{1/2} - \mA_k^{1/2}, \mI - \mA_{k+1}^{-1/2}\mA_k^{1/2}>
  =
  \eta^2\<\mA_{k+1}^{1/2} - \mA_k^{1/2}, \mI>   + \<\mA_k, \mA_{k+1}^{-1/2} - \mA_k^{-1/2}>
  \\&\at{uses the fact that $\mO \preceq \mA_k \preceq \mA_{k+1}$ and the L\"owner-Heinz theorem \citep[Theorem~2.6]{carlen2010trace}}\leq
  \eta^2\brs*{\normt{\mA_{k+1}^{1/2}} - \normt{\mA_k^{1/2}}}
  \at{uses the definition in \cref{eq:Psi}}=
  \eta \brs*{\Psi_{k+1}(m) - \Psi_k(m)}
  \at{uses \Cref{lem:Psi2}}\leq
  \eta^2\brs*{\normt{\sqrt{\mS_{k+1}}} - \normt{\sqrt{\mS_k}}},
\end{align*}
where \annotate.\qed

\proofsubsection{thm:FTRL}

We can upper-bound $\reg_K^+$ as follows:
\begin{align*}
  \reg_K^+
  &\at{uses \Cref{lem:regret+}}\leq
  \Psi_0^*(0) + \eta \distd{m_K} - \Psi_K^*(m_K)
  \\&
  + \msum_{k=1}^K \brs*{\Psi_k^*(m_{k-1}) - \Psi_{k-1}^*(m_{k-1})}
  - \msum_{k=0}^K \bg_{\Psi_k^*}(m_{k-1};m_k),
  \\&\at{uses the dominance property in \Cref{lem:Psi2}}\leq
  \Psi_0^*(0)
  + \msum_{k=1}^K \brs*{\Psi_k^*(m_{k-1}) - \Psi_{k-1}^*(m_{k-1})}
  - \msum_{k=0}^K \bg_{\Psi_k^*}(m_{k-1};m_k),
  \\&\at{uses the upper stability property in \Cref{lem:Psi2}}\leq
  \Psi_0^*(0)
  + \msum_{k=1}^K \eta\brs*{\normt{\sqrt{\mS_k}} - \normt{\sqrt{\mS_{k-1}}}}
  - \msum_{k=0}^K \bg_{\Psi_k^*}(m_{k-1};m_k),
  \\&\at{uses the smoothness property in \Cref{lem:Psi2}}\leq
  \Psi_0^*(0) + \eta\normt{\sqrt{\mS_K}}
  - \msum_{k=0}^K \tmfrac{1}{2\eta}\sqn{\nabla \Psi_k^*(m_{k-1}) - \nabla \Psi_k^*(m_k)}_{\mS_k^{1/2}}
  \\&\at{uses \cref{eq:FTRL2,eq:FTRL}, \Cref{lem:Psi}, and the assumption $x_0 = 0$}=
  \Psi_0^*(0) + \eta\normt{\sqrt{\mS_K}}
  - \tmfrac{1}{2\eta}\sqn{x_1 - x_0}_{\mS_0^{1/2}}
  - \msum_{k=1}^K \tmfrac{1}{2\eta}\sqn{\nabla \Psi_k^*(m_{k-1}) - x_{k+1}}_{\mS_k^{1/2}}
  \\&\at{uses \cref{eq:FTRL2}}\leq
  \Psi_0^*(0) + \eta\normt{\sqrt{\mS_K}}
  - \tmfrac{1}{2\eta}\sqn{x_1 - x_0}_{\mS_0^{1/2}}
  - \msum_{k=1}^K \tmfrac{1}{4\eta}\sqn{x_{k+1} - x_k}_{\mS_k^{1/2}}
  \\&
  + \msum_{k=1}^K \tmfrac{1}{2\eta}\sqn{\nabla \Psi_k^*(m_{k-1}) - \nabla \Psi_{k-1}^*(m_{k-1})}_{\mS_k^{1/2}}
  \\&\at{uses \Cref{lem:Psi3}}\leq
  \Psi_0^*(0) + \eta\normt{\sqrt{\mS_K}}
  - \msum_{k=0}^K \tmfrac{1}{4\eta}\sqn{x_{k+1} - x_k}_{\mS_k^{1/2}}
  + \msum_{k=1}^K \tfrac{1}{2}\eta\brs*{\normt{\sqrt{\mS_k}} - \normt{\sqrt{\mS_{k-1}}}}
  \\&\leq
  \Psi_0^*(0) + \tfrac{3}{2}\eta\normt{\sqrt{\mS_K}}
  - \msum_{k=0}^K \tmfrac{1}{4\eta}\sqn{x_{k+1} - x_k}_{\mS_k^{1/2}},
\end{align*}
where \annotate.\qed

\newpage

\proofsubsection{thm:online}

We can upper-bound $\reg_K$ as follows:
\begin{align*}
  \reg_K
  &=
  \reg_K^+ + \msum_{k=0}^K \<g_k,x_k - x_{k+1}>
  \\&\at{uses Young's inequality}\leq
  \reg_K^+ + \msum_{k=0}^K \brs*{\eta\sqn{g_k}_{\mS_k^{-1/2}} + \tmfrac{1}{4\eta}\sqn{x_{k+1} - x_k}_{\mS_k^{1/2}}}
  \\&\at{uses \Cref{thm:FTRL}}\leq
  \Psi_0^*(0) + \tfrac{3}{2}\eta\normt{\sqrt{\mS_K}}
  + \msum_{k=0}^K \eta\sqn{g_k}_{\mS_k^{-1/2}}
  \\&\at{uses \cref{line1:S,line1:init}, the definition in \cref{eq:Psi}, Lemma~3 of \citet{an2025asgo}, and the definition in \cref{eq:R}}\leq
  \eta\delta\dim{\cX} + \eta\distd{g_0} + \tfrac{3}{2}\eta\normt{\sqrt{\mS_K}}
  + \msum_{k=0}^K \eta\sqn{g_k}_{\mS_k^{-1/2}}
  \\&\at{uses the fact that $\mS_k^{-1/2} \in \cH$ due to \Cref{ass:H}, and the properties of the projection}=
  \eta\delta\dim{\cX} + \eta\distd{g_0} + \tfrac{3}{2}\eta\normt{\sqrt{\mS_K}}
  + \msum_{k=0}^K \eta\<\mS_k^{-1/2},\proj[\cH]{\out{g_k}}>
  \\&\at{uses \cref{line1:S}}=
  \eta\delta\dim{\cX} + \eta\distd{g_0} + \tfrac{3}{2}\eta\normt{\sqrt{\mS_K}}
  + \msum_{k=0}^K \eta\<\mS_k^{-1/2},\mS_k - \mS_{k-1}>
  \\&=
  \eta\delta\dim{\cX} + \eta\distd{g_0} + \tfrac{3}{2}\eta\normt{\sqrt{\mS_K}}
  \\&
  + \msum_{k=0}^K \eta\<\mS_k^{-1/2},(\mS_k^{1/2} - \mS_{k-1}^{1/2})(\mS_k^{1/2} + \mS_{k-1}^{1/2}) + \mS_{k-1}^{1/2}\mS_k^{1/2} - \mS_k^{1/2}\mS_{k-1}^{1/2}>
  \\&=
  \eta\delta\dim{\cX} + \eta\distd{g_0} + \tfrac{3}{2}\eta\normt{\sqrt{\mS_K}}
  + \msum_{k=0}^K \eta\<\mI + \mS_k^{-1/2}\mS_{k-1}^{1/2},\mS_k^{1/2} - \mS_{k-1}^{1/2}>
  \\&\at{uses H\"older's inequality for Schatten norms and the triangle inequality}\leq
  \eta\delta\dim{\cX} + \eta\distd{g_0} + \tfrac{3}{2}\eta\normt{\sqrt{\mS_K}}
  + \msum_{k=0}^K \eta\brs*{1 + \normo{\mS_k^{-1/2}\mS_{k-1}^{1/2}}}\normt{\mS_k^{1/2} - \mS_{k-1}^{1/2}}
  \\&\at{uses the fact that $\mS_k \succeq \mS_{k-1}$ due to \cref{line1:S} and \Cref{ass:H}}\leq
  \eta\delta\dim{\cX} + \eta\distd{g_0} + \tfrac{3}{2}\eta\normt{\sqrt{\mS_K}}
  + \msum_{k=0}^K 2\eta\brs*{\normt{\sqrt{\mS_k}} - \normt{\sqrt{\mS_{k-1}}}}
  \\&\leq
  \eta\delta\dim{\cX} + \eta\distd{g_0} + \tfrac{7}{2}\eta\normt{\sqrt{\mS_K}},
\end{align*}
where \annotate.\qed

\newpage

\section[Proofs for \crtCref{sec:convex}]{Proofs for \Cref{sec:convex}}

\proofsubsection{lem:function_gap}

First, we can lower-bound $\E{\<g_k,x_{k+1} - x>}$ for $k \in \N_0$ as follows:
\begin{align*}
  \mi{2}\E{\<g_k,x_{k+1} - x>}
  \\&=
  \E{\<g_k,x_{k+1} - x_k>} + \E{\<g_k,x_k - x>}
  \\&\at{uses \Cref{ass:grad2} and the fact that $x_k$ is independent of $\xi_k$}=
  \E{\<g_k,x_{k+1} - x_k>} + \E{\<\nabla f_k(x_k),x_k - x>}
  \\&\at{uses the convexity in \Cref{ass:f}}\geq
  \E{\<g_k,x_{k+1} - x_k>} + \E{f_k(x_k) - f_k(x)}
  \\&=
  \E{\<g_k,x_{k+1} - x_k>} + \E{f_k(x_{k+1}) + \<\nabla f_k(x_{k+1}),x_k - x_{k+1}> + \bg_{f_k}(x_k;x_{k+1}) - f_k(x)}
  \\&=
  \E{\<g_k - \nabla f_k(x_{k+1}),x_{k+1} - x_k>} + \E{\bg_{f_k}(x_k;x_{k+1}) + f_k(x_{k+1}) - f_k(x)}
  \\&\at{uses \cref{line2:g2}, \Cref{ass:grad2}, and the fact that $x_{k+1}$ and $x_k$ are independent of $\xi'_{k+1}$}=
  \E{\<g_k - \tg_{k+1},x_{k+1} - x_k>} + \E{\bg_{f_k}(x_k;x_{k+1}) + f_k(x_{k+1}) - f_k(x)}
  \\&\at{uses Young's inequality}\geq
  -\E*{\tfrac{1}{4}\eta^{-1}\sqn{x_{k+1} - x_k}_{\mS_{k+1}^{1/2}} + \eta\sqn{g_k - \tg_{k+1}}_{\mS_{k+1}^{-1/2}}}
  \\&
  + \E{\bg_{f_k}(x_k;x_{k+1}) + f_k(x_{k+1}) - f_k(x)}
\end{align*}
where \annotate. After rearranging and summing these inequalities for $k = 0, \dots, K$, we obtain the following:
\begin{align*}
  \mi{3}
  \msum_{k=0}^K\E{\bg_{f_k}(x_k;x_{k+1}) + f_k(x_{k+1}) - f_k(x)}
  \\&\at{uses the definition in \cref{eq:regret2}}\leq
  \E*{\msum_{k=0}^K \<g_k,x_{k+1} - x>}
  +\msum_{k=0}^K\E*{\tmfrac{1}{4\eta}\sqn{x_{k+1} - x_k}_{\mS_{k+1}^{1/2}} + \eta\sqn{g_k - \tg_{k+1}}_{\mS_{k+1}^{-1/2}}}
  \\&\at{is obtained by taking the maximum over $x \in Q_\cR$ and using the definition in \cref{eq:regret2}}\leq
  \E*{\reg_K^+}
  +\msum_{k=0}^K\E*{\tmfrac{1}{4\eta}\sqn{x_{k+1} - x_k}_{\mS_{k+1}^{1/2}} + \eta\sqn{g_k - \tg_{k+1}}_{\mS_{k+1}^{-1/2}}}
  \\&\at{uses \Cref{thm:FTRL}}\leq
  \E*{\Psi_0^*(0) + \tfrac{3}{2}\eta\normt{\sqrt{\mS_K}}}
  + \msum_{k=0}^K\E*{\tmfrac{1}{4\eta}\sqn{x_{k+1} - x_k}_{\mS_{k+1}^{1/2} - \mS_k^{1/2}} + \eta\sqn{g_k - \tg_{k+1}}_{\mS_{k+1}^{-1/2}}}
  \\&\at{uses the properties of the projection and the fact that $\mS_k^{-1/2}\in \cH$ and $\mS_{k+1}^{1/2} - \mS_k^{1/2}\in \cH$ due to \cref{line2:S,line2:init} and \Cref{ass:H}}=
  \E*{\Psi_0^*(0) + \tfrac{3}{2}\eta\normt{\sqrt{\mS_K}}}
  + \msum_{k=0}^K\E*{\tmfrac{1}{4\eta}\<\proj[\cH]{\out{x_{k+1} - x_k}},\mS_{k+1}^{1/2} - \mS_k^{1/2}>}
  \\&
  + \msum_{k=0}^K\E*{\eta\<\proj[\cH]{\out{g_k - \tg_{k+1}}},\mS_{k+1}^{-1/2}>}
  \\&\at{uses \cref{line2:S}, H\"older's inequality for Schatten norms, and the definition in \cref{eq:R}}\leq
  \E*{\Psi_0^*(0) + \tfrac{3}{2}\eta\normt{\sqrt{\mS_K}}}
  + \msum_{k=0}^K\E*{\tmfrac{1}{4\eta}\sqdist{x_{k+1} - x_k}\normt{\mS_{k+1}^{1/2} - \mS_k^{1/2}}}
  \\&
  + \msum_{k=0}^K\E*{\eta\<\mS_{k+1} - \mS_k,\mS_{k+1}^{-1/2}>}
  \\&\at{uses the triangle inequality, \cref{eq:FTRL2}, and the feasibility property in \Cref{lem:Psi2}}\leq
  \E*{\Psi_0^*(0) + \tfrac{3}{2}\eta\normt{\sqrt{\mS_K}}}
  + \msum_{k=0}^K\E*{\tfrac{1}{2}\eta\normt{\mS_{k+1}^{1/2} - \mS_k^{1/2}} + \eta\<\mS_{k+1} - \mS_k,\mS_{k+1}^{-1/2}>}
  \\&=
  \E*{\Psi_0^*(0) + \tfrac{3}{2}\eta\normt{\sqrt{\mS_K}}}
  + \msum_{k=0}^K\E*{\tfrac{1}{2}\eta\normt{\mS_{k+1}^{1/2} - \mS_k^{1/2}}}
  \\&
  + \msum_{k=0}^K\E*{\eta\<(\mS_{k+1}^{1/2} - \mS_k^{1/2})(\mS_{k+1}^{1/2} + \mS_k^{1/2}) + \mS_k^{1/2}\mS_{k+1}^{1/2} - \mS_{k+1}^{1/2}\mS_k^{1/2},\mS_{k+1}^{-1/2}>}
  \\&=
  \E*{\Psi_0^*(0) + \tfrac{3}{2}\eta\normt{\sqrt{\mS_K}}}
  + \msum_{k=0}^K\E*{\tfrac{1}{2}\eta\normt{\mS_{k+1}^{1/2} - \mS_k^{1/2}}}
  \\&
  + \msum_{k=0}^K\E*{\eta\<\mS_{k+1}^{1/2} - \mS_k^{1/2},\mI + \mS_{k+1}^{-1/2}\mS_k^{1/2}>}
  \\&\at{uses H\"older's inequality for Schatten norms and the triangle inequality}\leq
  \E*{\Psi_0^*(0) + \tfrac{3}{2}\eta\normt{\sqrt{\mS_K}}}
  + \msum_{k=0}^K\E*{\brs*{\tfrac{3}{2}\eta + \eta\normo{\mS_{k+1}^{-1/2}\mS_k^{1/2}}}\normt{\mS_{k+1}^{1/2} - \mS_k^{1/2}}}
  \\&\at{uses the fact that $\mS_{k+1} \succeq \mS_k$, which is implied by \cref{line2:S} and \Cref{ass:H}}\leq
  \E*{\Psi_0^*(0) + \tfrac{3}{2}\eta\normt{\sqrt{\mS_K}}}
  + \msum_{k=0}^K\E*{\tfrac{5}{2}\eta\normt{\mS_{k+1}^{1/2} - \mS_k^{1/2}}}
  \\&\at{uses the fact that $\mS_{k+1}^{1/2} \succeq \mS_k^{1/2}$ due to \cref{line2:S}, \Cref{ass:H}, and the L\"owner Heinz theorem \citep[Theorem~2.6]{carlen2010trace}}\leq
  \E*{\Psi_0^*(0) + 4\eta\normt{\sqrt{\mS_{K+1}}}}
  \\&\at{uses \cref{line2:init}}\leq
  \E*{\eta\delta \dim{\cX} + 4\eta\normt{\sqrt{\mS_{K+1}}}}
\end{align*}
where \annotate.\qed

\proofsubsection{lem:acceleration}

We can lower-bound $\msum_{k=0}^K\brs*{f_k(x_{k+1}) - f_k(x^*)}$ as follows:
\begin{align*}
  \msum_{k=0}^K\brs*{f_k(x_{k+1}) - f_k(x^*)}
  &\at{uses \cref{line2:ox} and the definition in \cref{eq:f_k}}=
  \msum_{k=0}^K\alpha_k^2\brs*{f(\ox_{k+1}) - f(x^*/\alpha_k + (1-1/\alpha_k)\ox_k)}
  \\&\at{uses \Cref{ass:f}}\geq
  \msum_{k=0}^K\alpha_k^2\brs*{f(\ox_{k+1}) - f(x^*)} - (\alpha_k^2 - \alpha_k)\brs{f(\ox_k) - f(x^*)}
  \\&\at{uses the fact that $f(\ox_k) - f(x^*) \geq 0$ which is implied by $\dist{\ox_k}\leq \eta$ and the assumption $\alpha_k^2 \leq \alpha_k + \alpha_{k-1}^2$}\geq
  \alpha_{K-1}^2 \brs*{f(\ox_K) - f(x^*)} - (\alpha_0^2 - \alpha_0)\brs*{f(\ox_0) - f(x^*)}
  \\&\at{uses the assumption $\alpha_0=1$}=
  \alpha_{K-1}^2 \brs*{f(\ox_K) - f(x^*)},
\end{align*}
where \annotate.\qed

\proofsubsection{lem:S}

We can upper-bound $\E*{\normt{\sqrt{\mS_{K+1}}}}$ as follows:
\begin{align*}
  \mi{3}\E*{\normt{\sqrt{\mS_{K+1}}}}
  \at{uses the fact that $\mS_{K+1} \in \Sp$}=
  \E*{\<\sqrt{\mS_{K+1}},\mI>}
  \\&\at{uses Lemma~3 of \citet{an2025asgo}}\leq
  \E*{\<\sqrt{\mS_0},\mI>} + \E*{\<\sqrt{\mS_{K+1} - \mS_0},\mI>}
  \\&\at{uses \cref{line2:init}}=
  \delta\dim{\cX} + \E*{\<\sqrt{\mS_{K+1} - \mS_0}\mB^{-1/2},\mB^{1/2}>}
  \\&\at{uses Young's inequality and the fact that $\normt{\mB} = 1$}\leq
  \delta\dim{\cX} + \E*{\sqrt{\<\mS_K - \mS_0, \mB^{-1}>}}
  \\&\at{uses \cref{line2:S,line2:init}}=
  \delta\dim{\cX} + \E*{\sqrt{\msum_{k=0}^K\<\proj[\cH]{\out{\tg_{k+1} - g_k}}, \mB^{-1}>}}
  \\&\at{uses the properties of the projection and the fact that $\mB^{-1} \in \cH$ due to \Cref{ass:H}}=
  \delta\dim{\cX} + \E*{\sqrt{\msum_{k=0}^{K}\<\out{\tg_{k+1} - g_k}, \mB^{-1}>}}
  \\&=
  \delta\dim{\cX} + \E*{\sqrt{\msum_{k=0}^K\sqn{\tg_{k+1} - g_k}_{\mB^{-1}}}}
  \\&\at{uses Jensen's inequality}\leq
  \delta\dim{\cX} + \sqrt{\msum_{k=0}^K\E*{\sqn{\tg_{k+1} - g_k}_{\mB^{-1}}}}
  \\&\at{uses \cref{line2:g1,line2:g2}, the definition in \cref{eq:f_k}, \Cref{ass:grad2}, and the fact that $\xi_k,\xi'_{k+1}$ are independent}\leq
  \delta\dim{\cX} +
  \sqrt{
    \msum_{k=0}^K 2\sigma^2\alpha_k^2
    +
    \E*{\msum_{k=0}^K \sqn{\nabla f_k(x_{k+1}) - \nabla f_k(x_k)}_{\mB^{-1}}
    }
  }
  \\&\at{uses \cref{eq:grad_diff}}\leq
  \delta\dim{\cX} +
  \sqrt{\msum_{k=0}^K 2\sigma^2\alpha_k^2
    +
    \E*{\msum_{k=0}^K
      \brs*{\tfrac{1+\nu}{\nu}}^{\frac{2\nu}{1+\nu}}
    \brs*{\cL\alpha_k^{1-\nu}}^{\frac{2}{1+\nu}}\brs*{\bg_{f_k}(x_k;x_{k+1})}^{\frac{2\nu}{1+\nu}}}
  }
  \\&\leq
  \delta\dim{\cX}
  +
  \sigma\sqrt{\msum_{k=0}^K 2\alpha_k^2}
  +
  \brs*{\tfrac{1+\nu}{\nu}}^{\frac{\nu}{1+\nu}}\cL^{\frac{1}{1+\nu}}\sqrt{\E*{\msum_{k=0}^K\brs*{\alpha_k^2}^{\frac{(1-\nu)}{1+\nu}}\brs*{\bg_{f_k}(x_k;x_{k+1})}^{\frac{2\nu}{1+\nu}}}}
  \\&\at{uses H\"older's inequality}\leq
  \delta\dim{\cX}
  +
  \sigma\sqrt{\msum_{k=0}^K 2\alpha_k^2}
  +
  \brs*{\tfrac{1+\nu}{\nu}}^{\frac{\nu}{1+\nu}}\cL^{\frac{1}{1+\nu}}\brs*{\msum_{k=0}^K \alpha_k^2}^{\frac{1-\nu}{2(1+\nu)}}\brs*{\E*{\msum_{k=0}^K \bg_{f_k}(x_k;x_{k+1})}}^{\frac{\nu}{1+\nu}}
  \\&\at{uses Young's inequality}\leq
  \delta\dim{\cX}
  +
  \sigma\brs*{\msum_{k=0}^K 2\alpha_k^2}^{\frac{1}{2}}
  +
  \E*{\tmfrac{1}{c}\msum_{k=0}^K \bg_{f_k}(x_k;x_{k+1})}
  +
  \tfrac{c^{\nu}}{1+\nu}\cL\eta^{\nu} \brs*{\msum_{k=0}^K \alpha_k^2}^{\frac{1-\nu}{2}}
  \\&\at{uses the assumptions $\alpha_k \leq \alpha_{k+1}$ and $c \geq 1$}\leq
  \delta\dim{\cX}
  +
  2\sigma\brs*{[K+1]\alpha_K^2}^{\frac{1}{2}}
  +
  \E*{\tmfrac{1}{c}\msum_{k=0}^{K-1}\bg_{f_k}(x_k;x_{k+1})}
  +
  c\cL\eta^{\nu} \brs*{[K+1]\alpha_K^2}^{\frac{1-\nu}{2}},
\end{align*}
where \annotate.\qed

\proofsubsection{thm:convex}

We combine \Cref{lem:function_gap} and \Cref{lem:S} with $c = 4$ and obtain the following inequality:
\begin{align*}
  \mi{10}\E*{\msum_{k=0}^K\bg_{f_k}(x_k;x_{k+1}) + f_k(x_{k+1}) - f_k(x^*)}
  \leq
  \E*{\msum_{k=0}^K \bg_{f_k}(x_k;x_{k+1})}
  \\&
  +
  5\cR\delta\dim{\cX}
  +
  8\sigma\cR\brs*{[K+1]\alpha_K^2}^{\frac{1}{2}}
  +
  16\cL\cR^{1+\nu} \brs*{[K+1]\alpha_K^2}^{\frac{1-\nu}{2}}.
\end{align*}
After rearranging, we get the following inequality
\begin{align*}
  \E*{\msum_{k=0}^K f_k(x_{k+1}) - f_k(x^*)}
  &\leq
  5\cR\delta\dim{\cX}
  +
  8\sigma\cR\brs*{[K+1]\alpha_K^2}^{\frac{1}{2}}
  +
  16\cL\cR^{1+\nu} \brs*{[K+1]\alpha_K^2}^{\frac{1-\nu}{2}}.
\end{align*}
Next, we use \Cref{lem:acceleration} and obtain the following inequality:
\begin{align*}
  \alpha_K^2\E*{f(\ox_{K+1}) - f(x^*)}
  \leq
  5\cR\delta\dim{\cX}
  +
  8\sigma\cR\brs*{[K+1]\alpha_K^2}^{\frac{1}{2}}
  +
  16\cL\cR^{1+\nu} \brs*{[K+1]\alpha_K^2}^{\frac{1-\nu}{2}}.
\end{align*}
It remains to use the assumption $\alpha_K = \frac{1}{2}(K + 2)$.\qed

\end{document}